\documentclass[12pt]{article}
\usepackage{amssymb}
\usepackage{amsmath,amsthm}
\usepackage[latin1]{inputenc}
\usepackage{hyperref}
\usepackage[thinlines]{easytable}
\usepackage{color}
\usepackage{enumerate}
\usepackage{graphicx}
 \newtheorem{remark}{Remark}

 \newtheorem{theorem}[remark]{Theorem}
 \newtheorem{proposition}[remark]{Proposition}
 \newtheorem{corollary}[remark]{Corollary}
 
 \newtheorem{example}[remark]{Example}
 \newtheorem{conjecture}[remark]{Conjecture}

\title{Defensive alliances in graphs: a survey}

\author{ Ismael Gonz\'alez Yero$^{1}$ and Juan A.
Rodr\'{\i}guez-Vel\'{a}zquez$^{2}$\\
    \\
$^1${\small Departamento de Matem\'aticas, Escuela Polit\'ecnica Superior de Algeciras}\\
{\small Universidad de C\'adiz,} {\small
Av. Ram\'on Puyol s/n, 11202 Algeciras, Spain.} \\ {\small
ismael.gonzalez\@@uca.es}\\
$^2${\small Departament d'Enginyeria Inform\`atica i Matem\`atiques}\\
{\small Universitat Rovira i Virgili,}  {\small Av. Pa\"{\i}sos
Catalans 26, 43007 Tarragona, Spain.} \\{\small
juanalberto.rodriguez\@@urv.cat}
\\
}

\begin{document}

\maketitle

\begin{abstract}
A set $S$ of vertices of a graph $G$ is a defensive $k$-alliance in $G$ if every vertex of $S$ has at least $k$ more neighbors inside of $S$ than outside. This is primarily an expository article surveying the principal known results on defensive alliances in graph.  Its seven sections are: Introduction, Computational complexity and realizability,  Defensive $k$-alliance number, Boundary defensive $k$-alliances, Defensive alliances in Cartesian product graphs,  Partitioning a graph  into defensive $k$-alliances, and Defensive $k$-alliance free sets.
\end{abstract}



\section{Introduction}

Alliances occur in a natural way in real life. General speaking, an alliance can be understood as a collection of elements sharing similar objectives or having similar properties among all elements of the collection. In this sense, there exist alliances like the following ones: a group of people united by a common friendship, or perhaps by a common goal; a group of plants belonging to the same botanical family; a group of companies sharing the same economic interest; a group of Twitter users following or being followed among themselves; a group of Facebook users sharing a common activity.

Alliances in graphs were described first  by Kristiansen {\em et al.} in \cite{alliancesOne}, where  alliances were classified into defensive, offensive or powerful. Defensive alliances in graphs were defined as a set of vertices of the graph such that every vertex of the alliance has at most one neighbor outside of the alliance than inside of the alliance. After this seminal paper,  the issue has been studied intensively. Remarkable examples are the articles \cite{kdaf1,kdaf5}, where the authors generalized  the concept of defensive alliance to defensive $k$-alliance as a set $S$ of vertices of a graph $G$ with the property that every vertex in $S$ has at least $k$ more neighbors in $S$ than it has outside of $S$.

Throughout this survey $G=(V,E)$ represents a undirected finite graph without loops and multiple edges with set of vertices $V$ and set of edges $E$. The
order of $G$ is $|V|=n(G)$ and the size $|E|=m(G)$ (If there is no ambiguity we will use only $n$ and $m$). We denote two adjacent vertices $u,v\in V$ by $u\sim v$ and in this case we say that $uv$ is an edge of $G$ or $uv\in E$. For a nonempty set $X\subseteq V$, and a vertex $v\in V$, $N_X(v)$ denotes the set of neighbors that $v$ has in $X$: $N_X(v):=\{u\in X: u\sim v\}$ and the degree of $v$ in $X$ is denoted by $\delta_{X}(v)=|N_{X}(v)|.$ In the case $X=V$ we will use only $N(v)$, which is also called the open neighborhood of a vertex $v\in V$, and $\delta(v)$ to denote the degree of $v$ in $G$. The close neighborhood of a vertex $v\in V$ is $N[v]=N(v)\cup \{v\}$. The minimum and maximum degree of $G$ are denoted by $\delta$ and $\Delta$, respectively.

Given  $k\in \{-\Delta,\dots,\Delta\}$, a nonempty set $S\subseteq V$ is a \emph{defensive $k$-alliance} in $G=(V,E)$ if %for every $ v\in S$,
\begin{equation}\label{cond-A-Defensiva}
\delta _S(v)\ge \delta_{\overline{S}}(v)+k, \quad \forall v\in S.
\end{equation}
Notice that equation (\ref{cond-A-Defensiva}) is equivalent to
\begin{equation}\label{cond-defensiva-grado2}
\delta (v)\ge 2\delta_{\overline{S}}(v)+k, \quad \forall v\in S.
\end{equation}
The minimum cardinality of a defensive $k$-alliance in $G$ is the {\em defensive $k$-alliance number} and it is denoted by $a_k(G)$. The case $k=-1$ corresponds to the standard \textit{defensive alliances} defined in \cite{alliancesOne}.
A set $S \subseteq V$ is a \textit{dominating set} in $G$ if for every vertex $v\in \overline{S}$,  $\delta_S(v)> 0$ (every
vertex in $\overline{S}$ is adjacent to at least one vertex in $S$). The \textit{domination number} of $G$, denoted by $\gamma(G)$, is the minimum cardinality of a dominating set in $G$ \cite{bookdom1}. A defensive $k$-alliance $S$ is called \emph{global} if it forms a dominating set. The minimum cardinality of a global defensive $k$-alliance in $G$ is the {\em global defensive $k$-alliance number} and it is denoted by $\gamma_k^d(G)$.

As a particular case of defensive alliance, in \cite{boundary-def} was defined and studied the limit case of equation
(\ref{cond-A-Defensiva}). In this sense, they defined a set $S\subset V$ as a {\em boundary defensive $k$-alliance} in $G$, $k\in \{-\Delta,\dots,\Delta\}$, if
\begin{equation}\label{cond-boundary-def}
\delta_S(v)=\delta_{\overline{S}}(v)+k, \quad \forall v\in S.
\end{equation}
A boundary defensive  $k$-alliance in $G$ is called {\em global} if it forms a dominating set in $G$. Notice that equation (\ref{cond-boundary-def}) is equivalent to
\begin{equation}\label{cond-boundary-def-degree}
\delta (v)= 2\delta_{S}(v)-k \quad, \forall v\in S.
\end{equation}
Note that there are graphs which does not contain any boundary defensive $k$-alliance for some values of $k$. For instance, the hypercube graph $Q_3$ has no boundary defensive $0$-alliances.

Defensive alliances have been studied in different ways. The first results about defensive alliances were presented in \cite{note,alliancesOne} and after that several results have been appearing in the literature, like those in \cite{G-araujo,bouzefrane-glob-def-tree,eroh-bounds,matamala,chellali,partitionTrees,partitionnumber,Favaron-ind-dom,ararat-ctw,GlobalalliancesOne,partitionGrid,
global-defe-star,spectral,planar,GArs,yrs,tesisiga,SBF,line}.
The complexity of computing minimum cardinality of defensive $k$-alliances in graphs was studied in \cite{complej1,complej2,jamieson-tesis,jamieson,SBF}, where it was proved that this is an NP-complete problem. A spectral study of alliances in graphs was presented in \cite{spectral,yrs}, where the authors obtained some bounds for the defensive alliance number in terms of the algebraic
connectivity, the Laplacian spectral radius and the spectral radius\footnote{The second smallest eigenvalue of the Laplacian matrix of a graph $G$ is called the algebraic connectivity of $G$. The largest eigenvalue of the adjacency matrix of $G$ is the spectral radius of $G$.} of the graph. The global defensive
alliances in trees and planar graphs were studied in \cite{bouzefrane-glob-def-tree,ararat-ctw} and \cite{planar}, respectively. The defensive alliances in regular graphs and circulant graphs were studied in \cite{G-araujo}. Moreover, the alliances in complement graphs, line graphs and weighted graphs were
studied in \cite{SBF}, \cite{yrs,line} and \cite{jamieson-dean}, respectively. Some relations between the independence number and the defensive alliances number of a graph were obtained in \cite{chellali,Favaron-ind-dom}. Also, the partitions of a graph into defensive $k$-alliances were investigated in \cite{partitionTrees,partitionnumber,partitionGrid,yero2}. Next we survey the principal known results about defensive alliances.

\section{Computational complexity and realizability}

The complexity of computing the minimum cardinality of a defensive $k$-alliance was studied in \cite{matamala,jamieson-tesis,jamieson,SBF}. Consider the following  decision problem (for any fixed $k$).\\

\noindent$\begin{tabular}{|l|}
  \hline
  \mbox{DEFENSIVE $k$-ALLIANCE PROBLEM}\\
  \mbox{INSTANCE: A graph $G = (V, E)$ and a positive integer $\ell < |V|$.}\\
  \mbox{PROBLEM:  Does $G$ have a $k$-defensive alliance of size at most $\ell$?}\\
  \hline
\end{tabular}$\\

\begin{theorem}{\em \cite{SBF}}
For any $k\in \{-\Delta,\dots,\Delta\}$, DEFENSIVE $k$-ALLIANCE PROBLEM is NP-complete.
\end{theorem}

The above result supplements and generalizes %(and at the same time, unifies)
known  results obtained in \cite{matamala,jamieson-tesis,jamieson} for $k=-1$.  Also, as shown in \cite{jamieson-tesis,jamieson}, DEFENSIVE $(-1)$-ALLIANCE PROBLEM is NP-complete,  even when restricted to split, chordal  or bipartite graphs.
\\

\noindent$\begin{tabular}{|l|}
  \hline
  \mbox{GLOBAL DEFENSIVE $k$-ALLIANCE PROBLEM}\\
  \mbox{INSTANCE: A graph $G = (V, E)$ and a positive integer $\ell < |V|$.}\\
  \mbox{PROBLEM:  Does $G$ have a global defensive alliance of size at most $\ell$?}\\
  \hline
\end{tabular}$
\\

Up to our knowledge, a general solution for this problem is still unknown and, as we can see below, for   $k=-1$ the problem is NP-complete.

\begin{theorem}{\em \cite{complej1,jamieson-tesis}} \label{GkDAProblem}
GLOBAL DEFENSIVE $(-1)$-ALLIANCE PROBLEM is NP-complete.
\end{theorem}

As shown in \cite{jamieson-tesis},  GLOBAL DEFENSIVE $(-1)$-ALLIANCE PROBLEM is NP-complete, even when restricted to  chordal graphs or bipartite graphs.

No we consider some realizability results. Since every global $(-1)$-alliance   is also a dominating set, we know that
$\gamma_{-1}^d(G)\ge \gamma(G)$ for any graph $G$. Every global $(-1)$-alliance   is also a defensive alliance, so
$\gamma_{-1}^d(G)\ge a_{-1}(G)$. In fact, as was shown in \cite{eroh-bounds}, any three positive integers satisfying these inequalities are
achievable as the $(-1)$-alliance number, the domination number, and the global $(-1)$-alliance number of some graph $G$.

\begin{theorem}{\em \cite{eroh-bounds}}
For any positive integers $a, b$ and $c$ with $a\le c$ and $b\le c$, there exists a connected graph G such that $a_{-1}(G) = a$, $\gamma(G) = b$ and $\gamma_{-1}^d(G) = c$.
\end{theorem}

Based simply on the definitions, the domination number, global $(-1)$-alliance number,
and global $0$-alliance number must satisfy
$\gamma(G)\le \gamma_{-1}^d(G)\le \gamma_{0}^d(G)$ for any graph
$G$. The following question was studied in \cite{eroh-bounds}: Given any three positive integers $a\le b\le c$, is there a graph $G$ so that
$\gamma(G)=a$, $ \gamma_{-1}^d(G)=b$ and $\gamma_{0}^d(G)=c$?

\begin{theorem}{\em \cite{eroh-bounds}}
Let $a, b$ and $c$ be three positive integers with $a\le b\le c$, $2\le b$ and $c\le\frac{ab+2b-a\left\lceil\frac{b}{a}\right\rceil}{2}$. Then there exists a graph $G$ such that $\gamma(G) = a$, $\gamma_{-1}^d(G) = b$ and $\gamma_0^d(G) = c$.
\end{theorem}

The next result concerns not only the minimum cardinality of a defensive $(-1)$-alliance, defensive $0$-alliance or a global defensive $(-1)$-alliance of a graph but also the subgraphs induced by these alliances.

\begin{theorem}{\em \cite{eroh-bounds}}
Given $1\le a\le b$ and any two connected graphs $H_1$ and $H_2$ with orders $a$ and $b$ respectively, there exists a connected graph $G$ with the following properties.
\begin{itemize}
\item $H_1$ is isomorphic to the subgraph induced by the only defensive alliance of $G$ that has minimum cardinality $a_{-1}(G)$.
\item $H_2$ is isomorphic to the subgraph induced by the only strong defensive alliance of $G$ that has minimum cardinality $a_0(G)$.
\end{itemize}
\end{theorem}

%As a consequence of the above theorem they obtained the following result.

\begin{corollary}{\em \cite{eroh-bounds}}
For any $1\le a\le b$, there exists a connected graph $G$ with $a= a_{-1}(G)\le b = a_0(G)$.
\end{corollary}

As the following result states, any connected graph is the subgraph induced by the unique
minimum global $(-1)$-alliance ($0$-alliance) of some graph.

\begin{theorem}{\em \cite{eroh-bounds}}
Given a connected graph $H$, there exists a connected graph $G$ for which $H$ is the subgraph induced by the unique global defensive $(-1)$-alliance $($respectively, $0$-alliance$)$ of $G$ with minimum cardinality  $\gamma^d_{-1}(G)$ $($respectively, $\gamma^d_{0}(G))$.
\end{theorem}

\section{Defensive $k$-alliance number}

According to the definitions, the domination number, global $k$-alliance number and alliance number must satisfy
\begin{equation}
\gamma_{k+1}^d(G)\ge \gamma_k^d(G)\ge \gamma(G)\quad
{\rm and } \quad \gamma_k^d(G)\ge a_k(G)\ge a_{k-1}(G)
\end{equation}
for any graph $G$. Now we present some results related to the monotony of $a_k(G)$ and $\gamma_k^d(G)$.

\begin{theorem}{\bf \cite{yrs}}\label{th3}
Let $G$ be a graph of minimum degree $\delta$ and maximum degree $\Delta$.
For every $k,r\in \mathbb{Z}$ such that $ -\delta \le k\le
\Delta$
 and $0 \le r \le \frac{k+\delta }{2}$,  $$  a_{_{k-2r}}(G)+r\leq a_{_k}(G).$$
\end{theorem}

The following two results are obtained directly from Theorem \ref{th3}.

\begin{corollary}{\bf \cite{yrs}}\label{coro1} Let $G$ be a graph of minimum degree $\delta$ and maximum degree $\Delta$ and let $t\in \mathbb{Z}$.
 \begin{itemize}
\item If $\frac{1-\delta }{2} \le t \le
\frac{\Delta-1}{2}$, then  $a_{2t-1}(G)+1\leq a_{2t+1}(G).$
%\quad {\rm  and} \quad
 \item  If $\frac{2-\delta}{2} \le t \le
\frac{\Delta}{2}$, then $ a_{2(t-1)}(G)+1 \leq a_{2t}(G).$
\end{itemize}
\end{corollary}

\begin{corollary}{\bf \cite{yrs}}\label{coro2}
Let $G$ be a graph of minimum degree $\delta$.
For every $k\in \{0,\dots,\delta \}$,
\begin{itemize}
\item if $k$ is even, then $a_{-k}(G)+\frac{k}{2}\le a_{0}(G) \leq a_{k}(G)-\frac{k}{2},$
\item if $k$ is odd, then $a_{-k}(G)+\frac{k-1}{2} \le a_{-1}(G)\leq a_{k}(G)-\frac{k+1}{2}.$
\end{itemize}
\end{corollary}

\begin{theorem}{\rm \cite{GArs}}\label{th1}
Let $S$ be a global defensive $k$-alliance of minimum cardinality in
$G$. If $W\subset S$ is a dominating set in $G$, then for
every $r\in \mathbb{Z}$ such that $0 \le r \le
\gamma_k^d(G)-|W|$,
$$  \gamma_{_{k-2r}}^d(G)+r\leq \gamma_{_k}^d(G).$$
\end{theorem}

The first bounds on the defensive alliance number appeared in \cite{note,alliancesOne}. For instance, the following results were obtained.

\begin{theorem}{\em \cite{note,alliancesOne}}
For any graph $G$ of order $n$ and minimum degree $\delta$,
$$a_{-1}^d(G)\le \min\left\{n-\left\lceil\frac{\delta}{2}\right\rceil,\left\lceil\frac{n}{2}\right\rceil\right\},$$
and also
$$a_{0}^d(G)\le \min\left\{ n-\left\lfloor\frac{\delta}{2}\right\rfloor ,\left\lfloor\frac{n}{2}\right\rfloor +1\right\}.$$
\end{theorem}
After that, some generalizations of the above results for the case of defensive $k$-alliances were presented in \cite{yrs}.

\begin{theorem}{\em \cite{yrs}}\label{n-delta-k-lower-upper}
Let $G$ be a graph of order $n$, maximum degree $\Delta$ and minimum degree $\delta$.
\begin{itemize}
\item For any $k\in \{-\delta,\dots,\Delta\}$,
$$\left\lceil\displaystyle\frac{\delta+k+2}{2}\right\rceil \le
a_k^d(G)\leq
n-\left\lfloor\displaystyle\frac{\delta-k}{2}\right\rfloor.$$

\item For any $k\in \{-\delta,\dots ,0\}$, $$a_k^d(G)\le \left\lceil
\displaystyle\frac{n+k+1}{2}\right\rceil.$$
\end{itemize}
\end{theorem}

The above bounds are achieved, for instance, for the complete graph $G=K_n$ for every $k\in \{1-n, \dots,  n-1\}$.

The global defensive $k$-alliance number has been also bounded using some basic parameters of the graphs like minimum and maximum degrees, size, etc.  For instance, it was shown in  \cite{GlobalalliancesOne} that
for any graph $G$ of order $n$ and minimum degree $\delta$,
$$
 \frac{\sqrt{4n+1}-1}{2}\le \gamma_{-1}^{d}(G)\le
n-\left\lceil\frac{d}{2}\right\rceil
$$
and
$$
 \sqrt{n}\le \gamma_0^{d}(G)\le
n-\left\lfloor\frac{d}{2}\right\rfloor.
$$
The next result generalizes the previous bounds to the case of  global defensive $k$-alliances.

\begin{theorem}{\em \cite{GArs}}\label{pilageneral}
Let $G$ be a graph of order $n$,  maximum degree $\Delta$ and minimum degree $\delta$. For any $k\in \{-\Delta,...,\Delta\}$,
$$\displaystyle\frac{\sqrt{4n+k^2}+k}{2} \le \gamma^{d}_k(G) \le
n-\left\lfloor\displaystyle\frac{\delta-k}{2}\right\rfloor.$$
\end{theorem}
The upper bound is attained, for instance, for the complete graph
$G=K_n$ for every $k\in \{1-n, \dots,  n-1\}$.  The lower bound
is attained, for instance, for the 3-cube graph $G=Q_3$,  in
the following cases: $\gamma_{-3}^d(Q_3)=2$ and $
\gamma_1^d(Q_3)=\gamma_0^d(Q_3)=4$.

It was shown in  \cite{GlobalalliancesOne} that for any bipartite
graph $G$ of order $n$ and maximum degree $\Delta$,
$$\gamma^{d}_{-1}(G)\ge \left\lceil
\frac{2n}{\Delta+3}\right\rceil \quad {\rm and } \quad
\gamma^{d}_0(G)\ge \left\lceil \frac{2n}{\Delta+2}\right\rceil .$$

The generalization of these bounds to the case of global defensive $k$-alliances is shown in the following theorem.

\begin{theorem}{\em \cite{GArs}}\label{cota-inf-def-global-kalliance}
For any graph $G$ of order $n$ and maximum degree $\Delta$ and for any $k\in \{-\Delta,...,\Delta\}$,
$$\gamma_k^d(G)\ge \left\lceil\frac{n}{\left\lfloor\frac{\Delta-k}{2}\right\rfloor+1}\right\rceil.$$
\end{theorem}

The above bound is tight. For instance, for the Petersen graph the
bound is attained for every $k$: $3\le\gamma_{-3}^d(G)$,
$4\le\gamma_{-2}^d(G)=\gamma_{-1}^d(G)$, $5\le
\gamma_0^d(G)=\gamma_1^d(G)$ and $10\le
\gamma_2^d(G)=\gamma_3^d(G)$. For the 3-cube graph
$G=Q_3$, the above theorem leads to the following exact values
of $\gamma_{k}^d(Q_3)$: $2\le\gamma_{-3}^d(Q_3)$, $4\le
\gamma_0^d(Q_3)=\gamma_1^d(Q_3)$ and $8\le \gamma_2^d(Q_3)=\gamma_3^d(Q_3)$.

\subsection{Defensive $k$-alliance number  of some particular graphs classes}

We begin this section with a resume of the values for the (global) defensive alliance number of some basic families of graphs. These results have been obtained in \cite{GlobalalliancesOne,alliancesOne}.\\

\begin{center}
\begin{TAB}(5ex,5ex,5ex,5ex,5ex,5ex,5ex){|c|c|c|c|c|}{|c|c|c|c|c|c|c|}
  \hline
  Graph $\;\;G\;\;$ & $\;\;a_{-1}(G)\;\;$ & $\;\;a_{0}(G)\;\;$ & $\;\;\gamma_{-1}(G)\;\;$ & $\;\;\gamma_0(G)\;\;$ \\
  $K_n$ & $\left\lfloor\frac{n+1}{2}\right\rfloor$ & $\left\lceil\frac{n+1}{2}\right\rceil$ & $\left\lfloor\frac{n+1}{2}\right\rfloor$ &  $\left\lceil\frac{n+1}{2}\right\rceil$ \\
  $P_n\; (n\ge 3)\; n\not\cong 2 (4)$ & 1 & 2 & $\left\lfloor\frac{n}{2}\right\rfloor+\left\lceil\frac{n}{4}\right\rceil-\left\lfloor\frac{n}{4}\right\rfloor$ & $\left\lfloor\frac{n}{2}\right\rfloor+\left\lceil\frac{n}{4}\right\rceil-\left\lfloor\frac{n}{4}\right\rfloor$ \\
  $P_n\; (n\ge 3)\; n\cong 2 (4)$ & 1 & 2 & $\left\lfloor\frac{n}{2}\right\rfloor+\left\lceil\frac{n}{4}\right\rceil-\left\lfloor\frac{n}{4}\right\rfloor$-1 & $\left\lfloor\frac{n}{2}\right\rfloor+\left\lceil\frac{n}{4}\right\rceil-\left\lfloor\frac{n}{4}\right\rfloor$\\
  $C_n\; (n\ge 3)$ & 2 & 2 & $\left\lfloor\frac{n}{2}\right\rfloor+\left\lceil\frac{n}{4}\right\rceil-\left\lfloor\frac{n}{4}\right\rfloor$ & $\left\lfloor\frac{n}{2}\right\rfloor+\left\lceil\frac{n}{4}\right\rceil-\left\lfloor\frac{n}{4}\right\rfloor$ \\
  $K_{1,s}\; (r\ge 2)$ & 1 & 1 & $\left\lfloor\frac{s}{2}\right\rfloor+1$ & $\left\lceil\frac{s}{2}\right\rceil+1$ \\
  $K_{r,s}\; (r,s\ge 2)$ & $\left\lfloor\frac{r}{2}\right\rfloor+\left\lfloor\frac{s}{2}\right\rfloor$ & $\left\lceil\frac{r}{2}\right\rceil+\left\lceil\frac{s}{2}\right\rceil$ & $\left\lfloor\frac{r}{2}\right\rfloor+\left\lfloor\frac{s}{2}\right\rfloor$ & $\left\lceil\frac{r}{2}\right\rceil+\left\lceil\frac{s}{2}\right\rceil$\\
  \end{TAB}
\end{center}

Defensive alliances in regular graphs and circulant graphs were studied in \cite{G-araujo}.
In order to present some results from \cite{G-araujo} it is necessary to introduce some notation.
Given a graph $G=(V,E)$ and  a subset $S\subset V$, the subgraph induced by $S$ will be denoted by $\langle S\rangle$.  The $(k, \delta)$-induced alliance is the set of graphs $H$ of order $t$, minimum degree $\delta_H\ge \left\lfloor\frac{\delta}{2}\right\rfloor$, and maximum degree $\Delta_H \le \delta$, with no proper subgraph of minimum degree greater than $\left\lfloor\frac{\delta}{2}\right\rfloor$. This set is denoted by $\mathcal{S}_{(t,\delta)}$.

\begin{theorem}{\em \cite{G-araujo}}
If $G$ is a $\delta$-regular graph, then $S$ is a critical alliance\footnote{A critical alliance is an alliance such that it does not contain other alliance as a proper subset.} of $G$ of cardinality $t$ if and
only if $\langle S\rangle\in \mathcal{S}_{(t,\delta)}$.
\end{theorem}

The $(6)$-regular graphs $G$ satisfying that $a_{-1}^d(G)\in \{4,5,6,7\}$ were characterized in \cite{G-araujo}. For the case of circulant graphs the following results were obtained in \cite{G-araujo}.

\begin{theorem}{\em \cite{G-araujo}}
Let $G = CR(n,M)$ be a circulant graph with $|M|$ generators.
\begin{itemize}
\item[{\rm (i)}] If $\delta=2|M|$, then $|M|+1\le a_{-1}^d(G)\le 2^{|M|}$.
\item[{\rm (ii)}] If $\delta=2|M|-1$, then $|M|\le a_{-1}^d(G)\le 2^{|M|-1}$.
\end{itemize}
\end{theorem}

As a consequence, it was obtained in \cite{G-araujo} that for the case of $|M|=3$, it is satisfied that $4\le a_{-1}^d(G)\le 8$. Moreover, the authors of that article characterized the circulant graphs $G$ such that $a_{-1}^d(G)\in \{4,5,6,7\}$.

An other class of graphs in which have been studied its defensive alliances is the case of planar graphs. For instance, \cite{planar} was dedicated to study defensive alliances in planar graphs, where are some results like the following one.

\begin{theorem}{\em \cite{planar}}
Let $G$ be a planar graph of order $n$.
\begin{itemize}
\item[{\rm (i)}] If $n > 6$, then $\gamma_{-1}^d(G)\ge
\left\lceil\frac{n+12}{8}\right\rceil.$
\item[{\rm (ii)}] If $n > 6$ and $G$ is a
triangle-free graph, then $\gamma_{-1}^d(G)\ge
\left\lceil\frac{n+8}{6}\right\rceil.$
\item[{\rm (iii)}] If $n > 4$, then $\gamma_{0}^d(G)\ge
\left\lceil\frac{n+12}{7}\right\rceil.$
\item[{\rm (iv)}] If $n > 4$ and $G$ is a triangle-free graph, then $\gamma_{0}^d(G)\ge
\left\lceil\frac{n+8}{5}\right\rceil.$
\end{itemize}
\end{theorem}

The following result concerns the particular case of trees.

\begin{theorem}{\rm \cite{GArs}} For any tree $T$ of order $n$,
$\gamma_{k}^{a}(T)\ge
\left\lceil\displaystyle\frac{n+2}{3-k}\right\rceil.$
\end{theorem}

The above bound is attained for $k\in \{-4,-3,-2,0,1\}$ in the case
of $G=K_{1,4}$. As a particular case of above theorem  we can
derive the following lower bounds obtained in \cite{GlobalalliancesOne}.

%Global defensive alliances in trees have been also studied separately, for instance, \cite{bouzefrane-glob-def-tree,ararat-ctw,GlobalalliancesOne} are examples of that.

\begin{theorem}{\em \cite{GlobalalliancesOne}}
If $T$ is a tree of order $n$, then

\begin{enumerate}[{\rm (i)}]
\item $\frac{n + 2}{4}\le \gamma_{-1}^d(T)\le \frac{3n}{5}$,
\item $\frac{n+2}{3}\le \gamma_{0}^d(T)\le \frac{3n}{4}$,
\end{enumerate}
and these bounds are sharp.
\end{theorem}

Similar results were obtained in \cite{bouzefrane-glob-def-tree} by using also the leaves and support vertices of the tree, where the authors also characterized the families of graphs achieving equality in the bounds.

\begin{theorem}{\em \cite{bouzefrane-glob-def-tree}}
Let $T$ be a tree of order $n\ge 2$ with $l$ leaves and $s$ support vertices. Then
\begin{enumerate}[{\rm (i)}]
\item $\gamma_{-1}^d(T)\ge \frac{(3n - l - s + 4)}{8}$,
\item $\gamma_{0}^d(T)\ge \frac{(3n - l - s + 4)}{6}$.
\end{enumerate}
\end{theorem}

A $t$-ary tree is a rooted tree where each node has at most $t$ children. A complete $t$-ary tree is a $t$-ary tree in which all the leaves have the same depth and all the nodes except the leaves have $t$ children. We let $T_{t,d}$ be the complete $t$-ary tree with depth/height $d$. With the above notation we present the following results obtained in \cite{ararat-ctw}.

\begin{theorem}{\em \cite{ararat-ctw}}\label{T2d}
Let $n$ be the order of $T_{2,d}$. Then  for any $d$, $$\gamma_{-1}^d(T_{2,d})=
\left\lceil\frac{2n}{5}\right\rceil.$$
\end{theorem}

\begin{theorem}{\em \cite{ararat-ctw}}\label{T3d} Let $d$ be an integer greater than three,
\begin{itemize}
\item[{\rm (i)}] If $d$ is odd, then $\gamma_{-1}^d(T_{3,d})=
\left\lfloor\frac{19n}{36}\right\rfloor$.
\item[{\rm (ii)}] If $d$ is even, then $\gamma_{-1}^d(T_{3,d})= \left\lceil\frac{19n}{36}\right\rceil$.
\end{itemize}
\end{theorem}

\begin{theorem}{\em \cite{ararat-ctw}}
For $d\ge 2$ and $t\ge 2$,
$$t^{d-1}\left\lceil\frac{t-1}{2}\right\rceil+t^{d-1}+t^{d-2}\le \gamma_{-1}^d(T_{t,d})\le
t^{d-1}\left\lceil\frac{t-1}{2}\right\rceil+t^{d-1}+t^{d-2}+t^{d-3}.$$
\end{theorem}

An efficient algorithm to determine the global defensive alliance numbers of trees was proposed in \cite{Chang2012479}, where the authors gave
formulas to compute the global defensive alliance numbers of complete $r$-ary trees for $r = 2, 3, 4.$ Since Theorems \ref{T2d} and \ref{T3d} provide  formulas for $r=2,3$, here we include the formula for $r=4$.

\begin{theorem}{\em \cite{Chang2012479}}
$\gamma_{-1}^d(T_{4,1})=3$, $\gamma_{-1}^d(T_{4,2})=9$ and for all
$d\ge 3$,

$$\gamma_{-1}^d(T_{4,d})=
\left\{\begin{array}{cc}
               \displaystyle\frac{577\times 4^{d-1}+47}{255}, & \mbox{if $d \equiv 0 ({\rm mod} 4)$} \\
               \\
                             \displaystyle\frac{577\times 4^{d-1}+443}{255}, & \mbox{if $d\equiv 1({\rm mod} 4)$} \\
                             \\
                               \displaystyle\frac{577\times 4^{d-1}-13}{255}, & \mbox{if $d\equiv 2({\rm mod} 4)$}
                               \\
                               \\
                                 \displaystyle\frac{577\times 4^{d-1}-52}{255}, & \mbox{if $d\equiv 3 ({\rm mod}  4)$.}
                              \end{array}
\right. $$
\end{theorem}

Consider the family $\xi$ of trees $T$, where $T$ is a star of odd order or $T$ is the tree
obtained from $K_{1,2t_1}$, $K_{1,2t_2}$, ...,$K_{1,2t_s}$,  and $tP_4$ (the disjoint union of $t$ copies of $P_4$) by
adding $s+t-1$ edges between leaves of these stars and paths in such a way that the center
of each star $K_{1,2t_i}$ is adjacent to at least $1 + t_i$ leaves in $T$ and each leaf of every copy of
$P_4$ is incident to at least one new edge, where $t\ge  0,$  $s \ge 2$ and $t_i\ge 2$ for $i = 1, 2, . . . , s.$
Note that each support vertex of each tree in $\xi$ must be adjacent with at least $3$ leaves.

\begin{theorem}{\em \cite{Chen2013}}
Let $T$ be a tree of order $n \ge  3$  with $s$ support vertices. Then
$$\gamma^d_{-1}(T) \le \frac{n+s}{2},$$
with equality if and only if $T \in \xi$.
\end{theorem}

\subsection{Relations between the (global) defensive $k$-alliance number and other invariants}

It is well-known that the algebraic connectivity of a graph is probably the most important information contained in the Laplacian spectrum. This eigenvalue is related to several important graph invariants and it imposes reasonably good bounds on the values of several parameters of graphs which are very hard to compute. Now we present a result about defensive alliances, obtained in \cite{yrs}.

\begin{theorem}{\em \cite{yrs}}\label{spectral-defe}
For any connected graph $G$  and for every $k\in
\{-\delta,\dots,\Delta\}$,  $$ a_k^d(G)\ge
\left\lceil\frac{n(\mu+k+1)}{n+\mu}\right\rceil.$$
\end{theorem}

The cases $k=-1$ and $k=0$ in the above theorem were studied previously in \cite{spectral}. Other relations between defensive alliances and the eigenvalues of a graph appeared in \cite{tesisiga}, in this case related to the spectral radius.

\begin{theorem}{\em \cite{tesisiga}}\label{lambda-defe}
For every graph $G$ of order $n$ and spectral radius $\lambda$,
$$\gamma_{k}^d(G)\ge \left\lceil\frac{n}{\lambda -k+ 1}\right\rceil.$$
\end{theorem}

The particular cases of the above theorem $k=-1$ and $k=0$ were studied previously in \cite{spectral}.

Some relationships between the independence number (independent domination number) and the global defensive alliance number of a graph were investigated in \cite{chellali,Favaron-ind-dom}. For instance, the following results were obtained there.

\begin{theorem}{\em \cite{chellali}}\label{citar-indep-num-2}
For any tree $T$, $\gamma_{-1}^d(T)\le \beta_0(T)$, and this bound is sharp.
\end{theorem}

\begin{theorem}{\em \cite{chellali}}
If $T$ is a tree of order $n\ge  3$ with $s$ support vertices, then
\begin{itemize}
\item[{\rm (i)}] $\gamma_0^d(G)\le \frac{3\beta_0(T)-1}{2}$,
\item[{\rm (ii)}] $\gamma_0^d(G)\le \beta_0(T)+s-1$.
\end{itemize}
\end{theorem}

In order to present some results from \cite{Favaron-ind-dom} we introduce some notation defined in the mentioned article.

$\mathcal{F}_1$ is the family of graphs obtained from a clique $S$ isomorphic to $K_t$ by attaching $t = \delta_S(u)+1$
leaves at each vertex $u\in S$.

$\mathcal{F}_2$ is the family of bipartite graphs obtained from a balanced complete bipartite
graph $S$ isomorphic to $K_{t,t}$ by attaching $t + 1$ leaves at each vertex $u\in S$.

$\mathcal{F}_3$  is the family of trees obtained from a tree $S$ by attaching a set $L_u$ of $\delta_S(u) + 1$
leaves at each vertex $u\in S$.

\begin{theorem}{\em \cite{Favaron-ind-dom}}$\,$\label{citar-indep-num-22}
\begin{itemize}
\item[{\rm (i)}] Every graph $G$ satisfies $i(G)\le (\gamma_{-1}^d(G))^2-\gamma_{-1}^d(G)+1$ with equality if and only if $G\in \mathcal{F}_1$.
\item[{\rm (ii)}] Every bipartite graph $G$ satisfies $i(G)\le \frac{(\gamma_{-1}^d(G))^2}{4}+\gamma_{-1}^d(G)$ with equality if and only if $G\in \mathcal{F}_2$.
\item[{\rm (iii)}] Every tree $G$ satisfies $i(G)\le 2\gamma_{-1}^d(G)-1$ with equality if and only if $G\in \mathcal{F}_3$.
\end{itemize}
\end{theorem}

Similarly to the above result, some relationships between the independent domination number and the global defensive $0$-alliance number of a graph were obtained in \cite{Favaron-ind-dom}.

\subsection{Complement graph and line graph}

As special cases of graphs in which their defensive alliances have been investigated, we have the complement graph and the line graph.

\begin{theorem}{\em \cite{SBF}}
If $G$ is a graph of order $n$ with maximum degree $\Delta$, then
$$\left\lceil\frac{n-\Delta+k+1}{2}\right\rceil\le a_k^d(\overline{G})\le \left\lceil\frac{n+\Delta+k+1}{2}\right\rceil.$$
\end{theorem}

\begin{theorem}{\em \cite{SBF}}
Let $G$ be a graph of order $n$ such that  $\gamma(G)> 3$ and $k\in \{-\delta,..., 0\}$. If the minimum defensive $k$-alliance in
$G$ is not global, then
$$a_k^d(\overline{G})\le \left\{\begin{array}{cc}
                                \displaystyle\left\lfloor\frac{3n+k+5}{4}-\frac{\gamma(G)+\gamma(\overline{G})}{2}\right\rfloor, & \mbox{if $n+k$ is odd} \\
                                \\
                                \displaystyle\left\lfloor\frac{3n+k+6}{4}-\frac{\gamma(G)+\gamma(\overline{G})}{2}\right\rfloor, & \mbox{if $n+k$ is even.}
                              \end{array}
\right.$$
\end{theorem}

Hereafter, we denote by ${\cal L}(G)$ the \textit{line graph}
of a simple graph $G$. Some of the next results are a generalization, to defensive $k$-alliances, of the previous ones obtained in \cite{line} on defensive (-1)-alliances and defensive $0$-alliances.

\begin{theorem}{\em \cite{yrs}}
For any graph $G$ of maximum degree $\Delta$, and for every $k\in \{2(1-\Delta),..., 0\}$,
$$a_k^d(\mathcal{L}(G))\le \Delta+\left\lceil\frac{k}{2}\right\rceil.$$
\end{theorem}

\begin{theorem}{\em \cite{yrs}}
Let $G=(V,E)$ be a simple graph of maximum degree $\Delta$. Let $v\in V$ such that $\delta(v)=\Delta$, let $\delta_v=\max\{\delta(u)\;:\;u\sim v\}$ and let $\delta_*=\min\{\delta_v\;:\;\delta(v)=\Delta\}$. For every $k\in \{2-\delta_*-\Delta,...,\Delta-\delta_*\}$,
$$a_k^d(\mathcal{L}(G))\le \left\lceil\frac{\Delta+\delta_*+k}{2}\right\rceil.$$
Moreover, if $\delta_1\ge \delta_2\ge ...\ge \delta_n$ is the degree sequence of $G$, then for every $k\in \{2-\delta_1-\delta_2,...,\delta_1+\delta_2-2\}$,
$$a_k^d(\mathcal{L}(G))\ge \left\lceil\frac{\delta_n+\delta_{n-1}+k}{2}\right\rceil.$$
\end{theorem}

As a consequence of the above results, the following interesting result was obtained in \cite{yrs}.

\begin{corollary}{\em \cite{yrs}}
For any $\delta$-regular graph $G$, $\delta>0$, and for every $k\in \{2(1-\delta),..., 0\}$,
$$a_k^d(\mathcal{L}(G))= \delta+\left\lceil\frac{k}{2}\right\rceil.$$
\end{corollary}

The cases $k=-1$ and $k=0$ in the above results were studied previously in \cite{line}.

We recall that a graph $G=(V, E)$ is a
($\delta_1,\delta_2)$-semiregular bipartite graph if the set $V$ can
be partitioned into two disjoint subsets $V_1, V_2$ such that if
$u\sim v$ then $u\in V_1$ and $v\in V_2$ and also,
$\delta(v)=\delta_1$ for every $v\in V_1$ and $\delta(v)=\delta_2$
for every $v\in V_2$.

\begin{corollary}{\em \cite{yrs}}
For any $(\delta_1,\delta_2)$-semiregular bipartite graph $G$,
$\delta_1>\delta_2$, and for every $k\in
\{2-\delta_1-\delta_2,\dots,\delta_1-\delta_2\}$,
 $$a_k^d({\cal L}(G))=\left\lceil
\frac{\delta_1+\delta_2+k}{2} \right\rceil.$$
\end{corollary}

We should point out that from the results shown  in the other
sections of this article on $a_k(G)$, we can derive some new results
on $a_k({\cal L}(G))$.

\section{Boundary defensive $k$-alliances}

Several basic properties of boundary defensive alliances were presented in \cite{boundary-def}.

\begin{remark}{\em \cite{boundary-def}}
Let $G$ be a simple graph and let $k\in \{-\Delta,\dots,\Delta\}$. If for every $v\in V$, $\delta(v)-k$ is an odd number, then $G$ does
not contain any boundary defensive  $k$-alliance.
\end{remark}

\begin{remark}{\em \cite{boundary-def}}
If $S$ is a  defensive $k$-alliance in $G$ and $\bar{S}$ is a global offensive $(-k)$-alliance in $G$, then $S$ is a  boundary defensive
$k$-alliance in $G$.
\end{remark}

\begin{theorem}{\em \cite{boundary-def}}\label{th111}
Let $G=(V,E)$ be a graph and let $S\subset V$. Let $m(\langle S\rangle)$ be the size of $\langle S\rangle$ and let $c$ be the number of edges of $G$ with one endpoint in $S$  and the other endpoint outside of $S$. If $S$ is a boundary defensive $k$-alliance in $G$, then
\begin{itemize}
\item[{\rm (i)}] $m(\langle S\rangle)=\displaystyle\frac{c+|S|k}{2}.$
\item[{\rm (ii)}] If $G$ is
a $\delta$-regular graph, then   $m(\langle S\rangle)
=\displaystyle\frac{|S|(\delta+k)}{4}$ and
$c=\displaystyle\frac{|S|(\delta-k)}{2}$.
\end{itemize}
\end{theorem}

Notice that if $S$ is a boundary defensive $k$-alliance in a graph $G$, then $a_k^d(G)\le |S|$. So, lower bounds for defensive $k$-alliance number are also lower bounds for the cardinality of any boundary defensive $k$-alliance. Moreover, upper bounds for the cardinality of any boundary defensive $k$-alliance are upper bounds for the defensive $k$-alliance number. For instance, the lower bound shown in Theorem \ref{n-delta-k-lower-upper} leads to a lower bound for the cardinality of any boundary defensive $k$-alliance. In the next result we state an upper bound for the cardinality of any boundary defensive $k$-alliance, which is the same obtained in Theorem \ref{n-delta-k-lower-upper} for the defensive $k$-alliance number.

\begin{remark}{\em \cite{boundary-def}}\label{cotainfsup}
If $S$ is a boundary defensive $k$-alliance in a graph $G$, then
$$\left\lceil\displaystyle\frac{\delta+k+2}{2}\right\rceil \le |S| \le\left\lfloor\displaystyle\frac{2n-\delta+k}{2}\right\rfloor.$$
\end{remark}

As the following corollary shows, the above bounds are tight.

\begin{corollary}{\em \cite{boundary-def}}\label{coro-kn-ak}
The cardinality of every boundary defensive $k$-alliance $S$ in the complete graph of order $n$ is $|S|=\frac{n+k+1}{2}$.
\end{corollary}

As a consequence of the above corollary it is concluded that the complete graph $G=K_n$ has boundary defensive $k$-alliances if and only if  $n+k+1$ is even.

The boundary defensive alliances were also related with the (Laplacian) spectrum of the graph as we can see below. The following theorems show the relationship between the algebraic connectivity (and the Laplacian spectral radius) of a graph and the cardinality of its boundary defensive  $k$-alliances.

\begin{theorem}{\em\cite{boundary-def}}\label{teo-mu-kdef}
Let $G$ be a connected graph.  If  $S$ is a boundary defensive
$k$-alliance in $G$, then
$$\left\lceil\frac{n(\mu-\left\lfloor\frac{\Delta-k}{2}\right\rfloor)}{\mu}\right\rceil \le |S| \le\left\lfloor\frac{n(\mu_*-\left\lceil\frac{\delta-k}{2}\right\rceil)}{\mu_*}\right\rfloor.$$
\end{theorem}

If $G=K_n$, then $\mu=\mu_*=n$ and $\Delta=\delta=n-1$. Therefore, the above theorem leads to the same result as Corollary \ref{coro-kn-ak}.

\begin{theorem}{\em\cite{boundary-def}}
Let $G$ be a connected graph. If  $S$ is a  boundary defensive $k$-alliance in $G$, then
$$\left\lceil\frac{n(\mu+k+2)-\mu}{2n}\right\rceil \le |S| \le n-\left\lceil\frac{n(\mu-k)-\mu}{2n}\right\rceil.$$
\end{theorem}

Notice that in the case of the complete graph  $G=K_n$,  the above theorem leads to Corollary \ref{coro-kn-ak}.

Boundary defensive $k$-alliances were also studied for the case of planar subgraphs. The Euler formula states that for a connected planar graph of order $n$, size $m$ and $f$  faces, $n - m + f = 2$. As a direct consequence of Theorem \ref{th111} and the Euler formula it is obtained the following result.

\begin{corollary}{\em\cite{boundary-def}}\label{coro1}
Let $G=(V,E)$ be a graph and let $S\subset V$. Let  $c$ be the number of edges of $G$ with one endpoint in $S$  and the other endpoint outside of $S$. If $S$ is a boundary defensive $k$-alliance in $G$  such that $\langle S\rangle$ is planar connected with $f$ faces, then
\begin{itemize}
\item[{\rm (i)}] $|S|=\displaystyle\frac{c+4-2f}{2-k}, $ for $k\neq 2.$
\item[{\rm (ii)}] If $G$ is a $\delta$-regular graph, then   $|S|=\displaystyle\frac{4f-8}{\delta+k-4}$ and $c=\displaystyle\frac{2(\delta-k)(f-2)}{\delta+k-4}$, for $k\in \{5-\delta, ..., \delta\}$.
\end{itemize}
\end{corollary}

\begin{theorem}{\em\cite{boundary-def}}
Let $G$ be a graph and let $S$ be a boundary defensive $k$-alliance
in $G$ such that $\langle S\rangle$ is planar connected with $f$
faces; then
$$|S|\le\left\lfloor\frac{\sqrt{16-8f+(n+k-2)^2}+n+k-2}{2}\right\rfloor.$$
\end{theorem}

The above bound is tight. For instance, the bound is attained for the complete graph $G=K_5$ where any set of cardinality four forms a boundary defensive $2$-alliance and $\langle S\rangle\cong K_4$ is planar with $f=4$ faces.

\begin{theorem}{\em\cite{boundary-def}}\label{teo-planar-def-1}
Let $G$ be a graph and let $S$ be a boundary defensive $k$-alliance in $G$ such that $\langle S\rangle$ is planar connected with $f>2$ faces.
\begin{itemize}
\item[{\rm (i)}] If $k\in \{5-\Delta,\dots,\Delta\}$, then $|S|\ge
\displaystyle\left\lceil\frac{4f-8}{\Delta+k-4}\right\rceil,$

\item[{\rm (ii)}] If $k\in \{5-\delta,\dots,\Delta\}$, then $ |S|\le
\displaystyle\left\lfloor\frac{4f-8}{\delta+k-4}\right\rfloor.$
\end{itemize}
\end{theorem}

By Corollary \ref{coro1} the above bounds are tight.

\section{Defensive alliances in Cartesian product graphs}
We recall that the \emph{Cartesian product of two graphs} $G=(V_1,E_1)$ and
$H=(V_2,E_2)$ is the graph $G\Box H$, such that
$V(G\Box H)=V_1\times V_2$ and two vertices $(a,b),(c,d)$ are adjacent in  $G\Box H$ if and only if  either
\begin{itemize}
\item $a=c$ and $bd\in E_2$, or
\item $b=d$ and $ac\in E_1$.
\end{itemize}

The study of defensive alliances in Cartesian product graphs was initiated  in \cite{alliancesOne}, where the authors obtained the following result.

\begin{theorem}{\em \cite{alliancesOne}}\label{alianzas-one-cartesiano}
For any Cartesian product graph $G\Box H$,
\begin{itemize}
\item[{\rm (i)}]$a_{-1}^d(G\Box H)\le \min\{a_{-1}^d(G)a_0^d(H),a_{0}^d(G)a_{-1}^d(H)\}.$
\item[{\rm (ii)}]$a_{0}^d(G\Box H)\le a_{0}^d(G)a_0^d(H).$
\end{itemize}
\end{theorem}

Let the graphs $G=(V_1,E_1)$ and $H=(V_2,E_2)$ and let $S\subset V_1\times V_2$ be a set of vertices of $G\Box H$. Let $P_{G_i}(S)$ the projection of the set $S$ over $G_i$. Then for every $u\in P_{G}(S)$ and every $v\in P_{H}(S)$, it is defined $X_u=\{(x,v)\in S\;:\;x=u\}$ and $Y_v=\{(u,y)\in S\; :\; y=v\}$.

\begin{theorem}{\em \cite{yero2}}
If $S\subset V_1\times V_2$ is a defensive $k$-alliance in $G\Box H$, then for  every $u\in P_{G}(S)$ and for every $v\in P_{H}(S)$, $P_{H}(X_u)$ and $P_{G}(Y_v)$ are a defensive $(k-\Delta_1)$-alliance in $H$ and a defensive $(k-\Delta_2)$-alliance in $G$, where $\Delta_1$ and $\Delta_2$ are the maximum degrees of $G$ and $H$, respectively.
\end{theorem}

Notice that $P_{H}(S)=\displaystyle\bigcup_{u\in P_{G}(S)}P_{H}(X_u)\; \mbox{\rm and} \;P_{G}(S)=\displaystyle\bigcup_{v\in P_{H}(S)}P_{G}(Y_v).$

 Also, as the union of defensive $k$-alliances in a graph is a defensive $k$-alliance in the graph, it is obtained the following consequence of the above result.

\begin{corollary}{\em \cite{yero2}}\label{teo-projections-def}
Let the graphs $G=(V_1,E_1)$ and $H=(V_2,E_2)$ of maximum degree $\Delta_1$ and $\Delta_2$, respectively. If $S\subset V_1\times V_2$ is a defensive $k$-alliance in $G\Box H$, then the projections $P_{G}(S)$ and $P_{H}(S)$ of $S$ over the graphs $G$ and $H$ are a defensive $(k-\Delta_2)$-alliance and a defensive $(k-\Delta_1)$-alliance in $G$ and $H$, respectively.
\end{corollary}

\begin{corollary}{\em \cite{yero2}}
Let the graphs $G=(V_1,E_1)$ and $H=(V_2,E_2)$ of maximum degree $\Delta_1$ and $\Delta_2$, respectively. If $G\Box H$ contains defensive $k$-alliances, then $G_i$ contains defensive $(k-\Delta_j)$-alliances, with $i,j\in \{1,2\}$, $i\ne j$ and, as a consequence,
$$a_k^d(G\Box H)\ge \max\{a_{k-\Delta_2}^d(G),a_{k-\Delta_1}^d(H)\}.$$
\end{corollary}

Also, in \cite{yero2} some relationships between $a_{k_1+k_2}^{d}(G\Box H)$ and $a_{k_i}^{d} (G_i)$, $i\in \{1,2\}$, were studied .

\begin{theorem}{\em \cite{yero2}}\label{CartesianNG-alliance}
For any graph $G$ and $H$, if $S_1$ is a defensive $k_1$-alliance in $G$ and $S_2$ is a defensive $k_2$-alliance in $H$, then $S_1\times S_2$ is a defensive $(k_1+k_2)$-alliance in $G\Box H$ and
$$a_{k_1+k_2}^{d}(G\Box H)\le a_{k_1}^{d} (G)a_{k_2}^{d}(H).$$
\end{theorem}

The bound of the above theorem is a general case of the results obtained in Theorem \ref{alianzas-one-cartesiano}.  Another interesting consequence of Theorem \ref{CartesianNG-alliance} is the following.

\begin{corollary}{\em \cite{yero2}}
Let $G$ and $H$ be two graphs of order $n_1$ and $n_2$ and  maximum  degree $\Delta_1$ and $\Delta_2$, respectively. Let $s\in
 \mathbb{Z}$ such that $\max\{\Delta_1,\Delta_2\}\le s\le \Delta_1+\Delta_2+k$. Then
$$a_{_{k-s}}^{d}(G\Box H)\le \min \{a_{k}^{d} (G),a_{k}^{d}(H)\}.$$
\end{corollary}

As a consequence of Theorem \ref{CartesianNG-alliance} it is obtained the following relationship between global defensive alliances in Cartesian product graphs and global defensive alliances in its factors.

\begin{corollary}{\em \cite{yero2}}\label{ThGanmaN2-alliance}
Let the graphs $G=(V_1,E_1)$, $H=(V_2,E_2)$ of minimum degree $\delta_1$, $\delta_2$ and maximum degrees $\Delta_1$ and $\Delta_2$, respectively.
\begin{itemize}
\item[{\rm (i)}]If  $G$ contains a global defensive $k_1$-alliance, then for every integer $k_2\in \{-\Delta_2,...,\delta_2\}$, $G\Box H$ contains a global defensive $(k_1+k_2)$-alliance and
    $$\gamma_{_{k_1+k_2}}^{d}(G\Box H)\le \gamma_{_{k_1}}^d(G)n_2.$$
\item[{\rm (ii)}]If  $H$ contains a global defensive $k_2$-alliance, then for every integer $k_1\in \{-\Delta_1,...,\delta_1\}$, $G\Box H$ contains a global defensive $(k_1+k_2)$-alliance and
$$\gamma_{_{k_1+k_2}}^{d}(G\Box H)\le \gamma_{_{k_2}}^d(H)n_1.$$
\end{itemize}
\end{corollary}

For a particular study of global defensive $(-1)$-alliances of  Cartesian product of paths and cycles we cite \cite{Chang2012479}.

\section{Partitioning a graph  into defensive $k$-alliances}

Other point of interest in investigating defensive alliances is related to graph partitions in which each set is formed by a defensive alliance. The partitions of a graph into defensive $(-1)$-alliances  were studied in \cite{partitionTrees,partitionnumber}. In these articles the concept of (global) defensive alliance partition number, ($\psi_{-1}^{gd}(G)$) $\psi_{-1}^d(G)$, was defined as the maximum number of sets in a partition of a graph such that every set of the partition is a (global) defensive $(-1)$-alliance.

\begin{theorem}{\em \cite{partitionnumber}}
Let $G$ be a connected graph of order $n\ge 3$. Then
$$1\le \psi_{-1}^{d}(G)\le \left\lfloor n+\frac{3}{2}-\frac{\sqrt{1+4n}}{2}\right\rfloor.$$
\end{theorem}

\begin{theorem}{\em \cite{partitionnumber}}
Let $G$ be a graph with minimum degree $\delta$. Then
$$\psi_{-1}^{d}(G)\le \left\lfloor\frac{n}{\left\lceil\frac{\delta+1}{2}\right\rceil}\right\rfloor.$$
\end{theorem}

Moreover, the partitions of  trees and grid graphs into (global) defensive $(-1)$-alliances, were studied in \cite{partitionTrees} and \cite{partitionGrid}, respectively.

\begin{theorem}{\em \cite{partitionTrees}}
Let $G$ be a connected graph with minimum degree $\delta$. Then
$$\psi_{-1}^{gd}(G)\le 1+\left\lfloor\frac{\delta}{2}\right\rfloor.$$
\end{theorem}

As a consequence of the above result, the following interesting result was obtained in \cite{partitionTrees}.

\begin{corollary}{\em \cite{partitionTrees}}
Let $T$ be a tree of order $n\ge 3$. Then
$1\le \psi_{-1}^{gd}(T)\le 2.$
\end{corollary}

Moreover, some families of trees satisfying that $\psi_{-1}^{gd}(T)=1$ or $\psi_{-1}^{gd}(T)=2$ were characterized in \cite{partitionTrees}. The following results for the class of grid graphs $P_r\square  P_c$ are known from \cite{partitionGrid}.

\begin{theorem}{\em \cite{partitionGrid}}
For $4\le r\le c$,
$$\psi_{-1}^d(P_r\square P_c)=\left\lfloor\frac{r-2}{2}\right\rceil\left\lfloor\frac{c-2}{2}\right\rceil+r+c-2.$$
\end{theorem}

\begin{theorem}{\em \cite{partitionGrid}}
For $2\le r\le c$,
$\psi_{-1}^{gd}(P_r\square P_c)=2.$
\end{theorem}

For any graph $G=(V,E)$, in \cite{yero2} is defined the (\emph{global}) \emph{defensive $k$-alliance partition number} of $G$, denoted by ($\psi_{k}^{gd}(G)$) $\psi_k^{d}(G)$, to be the maximum number of sets in a partition of $V$ such that each set of the partition is a (global) defensive
$k$-alliance, where $k\in\{-\Delta,...,\delta\}$.

Extreme cases are $\psi_{-\Delta}^{d}(G)=n$, where each set composed of one vertex is a defensive ($-\Delta$)-alliance, and $\psi_{\delta}^{d}(G)=1$ for the case of a connected $\delta$-regular graph where the whole vertex set of $G$ is the only defensive $\delta$-alliance. A graph $G$ is \emph{partitionable} into (global) defensive $k$-alliances if ($\psi_{k}^{gd}(G)\ge 2$) $\psi_{k}^{d}(G)\ge 2$. Hereafter we will say that ($\Pi_r^{gd}(G)$) $\Pi_r^d(G)$ is a partition of $G$ into $r$ (global) defensive $k$-alliances.

The following family of graphs was considered in \cite{yero2} to analyze the tightness of several of its results.

\begin{example}{\em \cite{yero2}}\label{Ej1} {\rm
Let $k$ and $r$ be integers such that $r>1$ and $r+k>0$  and let ${\cal H}$ be a family of graphs whose vertex set is $V=\cup_{i=1}^r{V_i}$ where, for every $V_i$, $\langle V_i\rangle\cong K_{r+k}$ and $\delta_{V_j}(v)=1$, for every $v\in V_i$ and $j\ne i$. Notice that $\{V_1,V_2,...,V_r\}$ is a partition
of the graphs belonging to ${\cal H}$ into $r$ global defensive $k$-alliances. A particular family of graphs included in ${\cal H}$ is $K_{r+k}\square K_r$.}
\end{example}

Hereafter, ${\cal H}$ will denote the family of graphs defined in the above example.

\begin{theorem}{\em \cite{yero2}}
For every graph $G$ partitionable into global defensive $k$-alliances,
\begin{itemize}
\item[{\rm (i)}] $\psi_{k}^{gd}(G) \le \lfloor\frac{\sqrt{k^2+4n}-k}{2}\rfloor, $
\item[{\rm (ii)}] $\psi_{k}^{gd}(G) \le \lfloor\frac{\delta-k+2}{2}\rfloor$.
\end{itemize}
\end{theorem}

The above bounds are attained, for instance, in the following cases: $\psi_{-1}^{gd}(K_4\square C_4)=4$, $\psi_{0}^{gd}(K_3\square C_4)=3$, $\psi_{1}^{gd}(K_2\square C_4)=2$  and  $\psi_1^{gd}(P)=2$, where $P$ denotes the Petersen graph.

\begin{remark}{\em \cite{yero2}}
For every $k\in \{1-\delta,...,\delta\}$, if $\psi_{k}^{gd}(G)\ge 2$, then
$$\gamma_{k}^{d}(G)+\psi_{k}^{gd}(G)\le \frac{n+4}{2}.$$
\end{remark}

Example of equality in the above relation is $\gamma_{-1}^{d}(C_4\square K_2)+\psi_{-1}^{gd}(C_4\square K_2) = 6.$

\begin{theorem}{\em \cite{yero2}}
If $\psi_k^{gd}(G)>2$, then, for every $l\in \{1,...,\psi_k^{gd}(G)-2\}$, there exists a subgraph, $G_l$, of $G$ of order $n(G_l)\le n(G)-l\gamma_k^{d}(G)$ such that $\displaystyle \psi_{l+k}^{gd}(G_l)+l\ge \psi_k^{gd}(G).$
\end{theorem}

One example where $\displaystyle \psi_{l+k}^{gd}(G_l)+l=\psi_k^{gd}(G)$ and $n(G_l)= n(G)-l\gamma_k^{d}(G)$ is the following. Let $G=K_4\square C_4$, the Cartesian product of the complete graph $K_4$ by the cycle graph $C_4$. $\psi_{-1}^{gd}(K_4\square C_4)=4$ and we can take each set of $\Pi_4^{gd}(K_4\square C_4)$ as the vertex set of a copy of $C_4$, so $G=K_3\square C_4$ and $H=K_2\square C_4$ (the 3-cube graph). Hence, $4=\psi_{-1}^{gd}(K_4\square C_4)=\psi_{0}^{gd}(K_3\square
C_4)+1=\psi_{1}^{gd}(K_2\square C_4)+2$ and $8=n(K_2\square C_4)=n(K_3\square C_4)-\gamma_{-1}^{d}(K_3\square C_4)=[n(K_4\square C_4)-\gamma_{-1}^{d}(K_4\square C_4)]-\gamma_{-1}^{d}(K_3\square C_4)=n(K_4\square C_4)-2\gamma_{-1}^{d}(K_4\square C_4)=16-2\cdot 4.$

\begin{theorem}\label{cortearistas}{\em \cite{yero2}}
Let $C_{(r,k)}^{gd} (G)$ be the minimum number of edges having its endpoints in different sets of a partition of $G$ into $r\ge 2$ global defensive $k$-alliances. Then
\begin{itemize}
\item[{\rm (i)}] $C_{(r,k)}^{gd} (G)\ge \frac{1}{2}r(r-1)\gamma_k^d(G)$,
\item[{\rm (ii)}]  $C_{(r,k)}^{gd} (G)\ge \frac{1}{2}r(r-1)(r+k)$,
\item[{\rm (iii)}]  $C_{(r,k)}^{gd} (G) \le \frac{2m-nk}{4}.$
\item[{\rm (iv)}] $C_{(r,k)}^{gd} (G) = \frac{1}{2}r(r-1)\gamma_k^d(G)=\frac{1}{2}r(r-1)(r+k)=\frac{2m-nk}{4}$ if and only if $G \in {\cal H}$.
\end{itemize}
\end{theorem}

From Theorem \ref{n-delta-k-lower-upper} is obtained that
\begin{equation}\label{cotainf}a_{k}^{d}(G)\ge \left\lceil
\frac{\delta+k+2}{2}\right\rceil.
\end{equation}

By Theorem \ref{cortearistas}  and equation (\ref{cotainf}) are obtained the following two necessary conditions for the existence of a partition of a graph into $r$ global defensive $k$-alliances.

\begin{corollary}{\em \cite{yero2}}
If for a graph $G$, $k>\frac{2m-r(r-1)(\delta+2)}{n+r(r-1)}$ or $k>\frac{2(m-r^2(r-1))}{n+2r(r-1)}$, then $G$ cannot be partitioned into $r$ global defensive $k$-alliances.
\end{corollary}

\subsection{Partitioning a graph into boundary defensive $k$-alliances}

Let $G=(V,E)$ be a graph and let $\Pi_r^d(G) =\{S_1,S_2,...S_r\}$ be a partition of $V$ into $r$  boundary defensive $k$-alliances. Suppose $x=\displaystyle\max_{1\le i\le r}\{|S_i|\}$ and $y=\displaystyle\min_{1\le i\le r}\{|S_i|\}$. Thus, $ \frac{n}{x}\le r \le \frac{n}{y}.$ Examples of bounds of $r$ are the following two corollaries.

As a consequence of Remark \ref{cotainfsup} the following bounds are obtained .

\begin{corollary}{\em \cite{yero2}}
If $G$ can be partitioned into $r$ boundary defensive $k$-alliances, then $$\displaystyle\frac{2n}{2n-\delta+k}\le r\le \frac{2n}{\delta+k+2}.$$
\end{corollary}

The above bounds are tight. For instance, if $n$ is even, each pair of vertices of $K_n$ forms a boundary defensive $(3-n)$-alliance. Thus, $K_n$ can be partitioned into $\frac{n}{2}$ of these alliances.

As a consequence of Theorem \ref{teo-mu-kdef} the following result is obtained .

\begin{corollary}{\em \cite{yero2}}\label{cota2-para r}
If $G$ can be partitioned into $r$ boundary defensive $k$-alliances, then
$$\displaystyle\frac{2\mu_*}{2\mu_*-\delta+k}\le r\le\frac{2\mu}{2\mu-\Delta+k}.$$
\end{corollary}

The above bounds are tight. By Corollary \ref{cota2-para r} it is concluded, for instance, that if the Petersen graph can be partitioned into $r$ boundary defensive $k$-alliances, then $k=1$ and $r=2$ (in this case $\Delta=\delta=3$, $\mu =2$ and $\mu_*=5$).

\begin{theorem}{\em \cite{yero2}}\label{partition-two-def}
Let $G=(V,E)$ be a graph and let $M\subset E$ be a cut set partitioning $V$ into two  boundary defensive $k$-alliances $S$ and $\overline{S}$, where $k\neq \Delta$ and $k\neq \delta$. Then
$$\displaystyle\left\lceil\frac{2m-kn}{2(\Delta-k)}\right\rceil\le |S|\le \left\lfloor\frac{2m-kn}{2(\delta-k)}\right\rfloor\; \mbox{\rm and}\;
|M|=\displaystyle\frac{2m-kn}{4}.$$
\end{theorem}

\begin{corollary}{\em \cite{yero2}}\label{coroPartitionDefensiveRegular}
Let $G=(V,E)$ be a $\delta$-regular graph and let $M\subset E$ be a cut set partitioning $V$ into two boundary defensive $k$-alliances
$S$ and $\overline{S}$. Then $|S|=\frac{n}{2}$ and  $|M|=\frac{n(\delta-k)}{4}$.
\end{corollary}

\begin{theorem}{\em \cite{yero2}}\label{partition-mu-def}
If $\{X,Y\}$ is a partition of $V$ into two boundary defensive $k$-alliances in $G=(V,E)$, then, without loss of generality, $$\left\lceil\sqrt{\frac{n(kn-2m+n\mu)}{4\mu}}+\frac{n}{2}\right\rceil\le |X|\le \left\lfloor\sqrt{\frac{n(kn-2m+n\mu_*)}{4\mu_*}}+\frac{n}{2}\right\rfloor$$
and
$$\left\lceil\frac{n}{2}-\sqrt{\frac{n(kn-2m+n\mu_*)}{4\mu_*}}\right\rceil\le |Y|\le\left\lfloor\frac{n}{2}-\sqrt{\frac{n(kn-2m+n\mu)}{4\mu}}\right\rfloor.$$
\end{theorem}

By Corollary \ref{coroPartitionDefensiveRegular} and Theorem \ref{partition-mu-def} it is obtained the following interesting consequence.

\begin{theorem}{\em \cite{yero2}}
Let $G=(V,E)$ be a $\delta$-regular graph. If $G$ is partitionable into two boundary defensive $k$-alliances, then the algebraic connectivity of $G$ is $\mu=\delta-k$ $($an even number$)$.
\end{theorem}

By the above necessary condition of existence of a partition of $V$ into two boundary defensive $k$-alliances it follows that,  for instance, the icosahedron cannot be partitioned into two boundary defensive $k$-alliances, because its algebraic connectivity is $\mu=5-\sqrt{5}\not\in \mathbb{Z}$. Moreover, the Petersen graph can only be partitioned into two boundary defensive $k$-alliances for the case of $k=1$, because $\delta=3$ and $\mu=2$.

\subsection{Partitioning $G\Box H$ into defensive $k$-alliances}

In this subsection some relationships between $\psi_{k_1+k_2}^{d}(G\Box H)$ and $\psi_{k_i}^{d}(G_i)$, $i\in \{1,2\}$, are presented. From Theorem \ref{CartesianNG-alliance} it follows that if $G$ contains a defensive $k_1$-alliance and $H$ contains a defensive $k_2$-alliance, then $G\Box H$ contains a defensive $(k_1+k_2)$-alliance. Therefore, the following result is obtained.

\begin{theorem}{\em \cite{yero2}}\label{CartesianNG}
For any graphs $G$ and $H$, if there exists a partition of $G_i$ into defensive $k_i$-alliances, $i\in \{1,2\}$, then there exists a partition of $G\Box H$ into defensive $(k_1+k_2)$-alliances and $$\psi_{k_1+k_2}^{d}(G\Box H) \ge \psi_{k_1}^{d}(G)\psi_{k_2}^{d}(H).$$
\end{theorem}

In the particular case of the Petersen graph, $P$, and the $3$-cube graph, $Q_3$, it follows $\psi_{-2}^{d}(P\square Q_3) =20= \psi_{-1}^{d}(P)\psi_{-1}^{d}(Q_3)$ and $5=\psi_{2}^{d}(P\square Q_3) > \psi_{1}^{d}(P)\psi_{1}^{d}(Q_3)=4$. Notice that  Theorem \ref{CartesianNG} leads to  $\psi_{2k}^{d}(G\Box H)\ge \psi_{k}^{d}(G)\psi_{k}^{d}(H).$ Another interesting consequence of Theorem \ref{CartesianNG} is the following.

\begin{corollary} {\em \cite{yero2}}
Let $G_i$ be a graph of order $n_i$ and  maximum  degree $\Delta_i$, $i\in \{1,2\}$. Let $s\in  \mathbb{Z}$ such that  $\max\{\Delta_1,\Delta_2\}\le s\le \Delta_1+\Delta_2+k$. Then $$\psi_{k-s}^{d}(G\Box H) \ge \max \{n_2\psi_k^{d}(G),n_1\psi_k^{d}(H)\}.$$
\end{corollary}

As example of equality we take $G=P$, $H=Q_3$, $k=1$ and $s=3$. In such a case, $20=\psi_{-2}^{d}(P\square Q_3) = \max \{8\psi_1^{d}(P),10\psi_1^{d}(Q_3)\}=\max\{16,20\}$.

At next we presented some results about global defensive $k$-alliances.

\begin{theorem}{\em \cite{yero2}}\label{ThGlobCartesian}
Let $\Pi_{r_i}^{gd}(G_i)$ be a partition of a graph $G_i$, of order $n_i$, into $r_i\ge 1$ global defensive $k_i$-alliances, $i\in \{1,2\}$, $r_1\le r_2$. Let $x_i=\displaystyle\min_{X\in \Pi_{r_i}^{gd}(G_i)}\{|X|\}$. Then,
\begin{itemize}
\item[{\rm (i)}] $\displaystyle\gamma_{_{k_1+k_2}}^{d}(G_1\Box G_2)\le \min\left\{x_1n_2,x_2n_1 \right\},$
\item[{\rm (ii)}] $\displaystyle\psi_{k_1+k_2}^{gd}(G_1\Box G_2) \ge \max\left\{\psi_{k_1}^{gd}(G_1),\psi_{k_2}^{gd}(G_2)\right\}$.
\end{itemize}
\end{theorem}

\begin{corollary}{\em \cite{yero2}}\label{coroProduct}
If  $G_i$ is a graph  of order $n_i$ such that $\psi_{k_i}^{gd}(G_i)\ge 1$, $i\in \{1,2\}$, then $$\displaystyle\gamma_{_{k_1+k_2}}^{d}(G_1\Box G_2)\le
\frac{n_1n_2}{\max_{i\in \{1,2\}}\left\{\psi_{k_i}^{gd}(G_i)\right\}}.$$
\end{corollary}

For the graph $C_4\square Q_3$, by taking $k_1=0$ and $k_2=1$, equalities in Theorem \ref{ThGlobCartesian} and Corollary \ref{coroProduct} are obtained.

\section{Defensive $k$-alliance free sets}

A set $Y \subseteq V$ is a \emph{defensive}  {\em $k$-alliance cover}, $k$-dac, if for all defensive
 $k$-alliance $S$, $S\cap Y\neq\emptyset$, \emph{i.e.}, $Y$ contains at least one vertex from each defensive
$k$-alliance of $G$. A $k$-dac  set $Y$ is \emph{minimal} if no proper subset of $Y$ is a defensive
$k$-alliance cover set. A \emph{minimum} $k$-dac  set is a minimal cover set of smallest cardinality. Also, a set $X\subseteq V$ is {\em defensive} $k$-{\em alliance free set},  $k$-daf, if for all defensive $k$-alliance $S$, $S\setminus X\neq\emptyset$, \emph{i.e.}, $X$ does not contain any defensive $k$-alliance as a subset. A $k$-daf  set $X$ is \emph{maximal} if it is not a proper subset of any defensive
$k$-alliance free set. A \emph{maximum} $k$-daf set is a maximal free set of biggest cardinality.

Hereafter, if there is no restriction on the values of k, we assume that  $k\in \{-\Delta,..., \Delta\}$.

\begin{theorem}{\em \cite{kdaf1,kdaf2}}\label{thDual}
\mbox{}
\begin{enumerate}[{\em (i)}]
\item  $X$ is a defensive  $k$-alliance cover set
if and only if $\overline{X}$ is defensive
$k$-alliance free set.

\item If $X$ is a minimal $k$-dac set then, for all $v\in X$,
there exists a defensive  $k$-alliance $S_v$ for
which $S_v \cap X=\{ v \}$.

\item  If X is a maximal k-daf set, then, for all $v\in
\overline{X}$, there exists $S_v\subseteq X$ such that $S_v \cup \
\{v\}$ is a defensive  $k$-alliance.
\end{enumerate}
\end{theorem}

Associated with the characteristic sets defined above we have the
following invariants:

\begin{itemize}

\item[] $\phi_k(G)$:  cardinality of a maximum $k$-daf set in $G$.

\item[] $\zeta_k(G)$: cardinality of a minimum $k$-dac  set in
$G$.

\end{itemize}

%%%%%%%%%%%%%%%%%%%%%%%%%%%%%%%%%%%%%
%%%%%%%%%%%%%%%%%%%%%%%%%%%%%%%%%%

The following corollary is a direct consequence of Theorem \ref{thDual} (i).

\begin{corollary}{\em \cite{kdaf1,kdaf2}}
$\phi_k(G)+\zeta_k(G)=n.$
\end{corollary}

Our next result leads to a property related to the monotony of $\phi_{k}(G)$.

\begin{theorem}{\em \cite{Rodriguez-Velazquez2011}}
If $X$ is a $k$-daf set and  $v\in \overline{X}$, then $X\cup \{v\}$
is  $(k+2)-daf$.
\end{theorem}

\begin{corollary}{\em \cite{Rodriguez-Velazquez2011}}
 For every $k\in \{-\Delta,...,\Delta-2\}$ and $r\in
\left\{1,...,\lfloor\frac{\Delta-k}{2}\rfloor\right\}$,
$$\phi_{k}(G)+r\le \phi_{k+2r}(G).$$
\end{corollary}

Now we point out some known bounds  on $\phi_{k}(G)$ and one conjecture related to one of these bounds.

\begin{theorem}{\em \cite{kdaf1,kdaf2}}
 For any connected graph  $G$ and $k\in \{0,...,\Delta\}$,

  $$\phi_k(G)\ge
\left\lfloor\frac{n}{2}\right\rfloor+\left\lfloor\frac{k}{2}\right\rfloor .$$
\end{theorem}
%The authors of \cite{kdaf1,kdaf2} proposed  the following possible extension of the above theorem.

\begin{conjecture}{\em \cite{kdaf1,kdaf2}}
 For any connected graph  $G$ and $-\delta \le k\le \Delta $,

  $$\phi_k(G)\ge
\left\lfloor\frac{n}{2}\right\rfloor+\left\lfloor\frac{k}{2}\right\rfloor .$$
\end{conjecture}

 The next
result shows other bounds on $\phi_k(G)$.

\begin{theorem} \label{cotaconnectivity}
 For any connected graph  $G$  and $-\Delta \le k \le \Delta$,
$$\displaystyle\left\lceil\frac{n(k+\mu)-\mu}{n+\mu}\right\rceil
\le\phi_{k}(G)\le  \left\lfloor\frac{2n+k-\delta-1}{2}
\right\rfloor,$$ where $\mu$ denotes the algebraic connectivity of
$G$.
\end{theorem}

The above bound is sharp as we can check, for instance, for the
complete graph $G=K_n$.  As the algebraic connectivity of $K_n$
is $\mu=n$, the above theorem gives the exact value of
$\phi_{k}(K_{n})= \left\lceil\frac{n+k-1}{2}\right\rceil.$

\begin{theorem}\label{thOf}
For any connected graph  $G$ and $-\Delta \le k \le \Delta$, $$\displaystyle\zeta_k(G)\le
\displaystyle\frac{n}{\mu_*} \left( \mu_*-\left\lceil\frac{\delta
+k}{2}\right\rceil\right),$$ where $\mu_*$ denotes the Laplacian
spectral radius of $G$.
\end{theorem}

The above bound is tight. For instance, we consider the
complete graph $G=K_n$ for which the Laplacian spectral radius  is $\mu_*=n$. In such a case, the above theorem gives the exact value
$\zeta_{k}(K_{n})= \left\lceil\frac{n-k}{2}\right\rceil.$

Now we state the following fact that will be useful for an easy understanding of some examples in the next subsection.

\begin{proposition}{\em \cite{Yero2013}}\label{remarktree}
Let $G$ be a graph of order $n$ and maximum degree $\Delta$. Then $\phi_k^d(G)=n$, for each of the following  cases:
\begin{itemize}
\item[{\rm (i)}] $G$ is a tree of maximum degree $\Delta\ge 2$ and $k\in \{2,...,\Delta\}$.
\item[{\rm (ii)}] $G$ is a planar graph of maximum degree $\Delta\ge 6$ and $k\in \{6,...,\Delta\}$.
\item[{\rm (iii)}] $G$ is a planar triangle-free graph of maximum degree $\Delta\ge 4$ and $k\in \{4,...,\Delta\}$.
\end{itemize}
\end{proposition}

\subsection{Defensive $k$-alliance free sets in Cartesian product graphs} \label{Defensive}

To begin with the study we present the following straightforward result.
 \begin{remark}{\em \cite{Yero2013}} \label{remark1}
 Let $G_i$ be a graph of order $n_i$, minimum degree $\delta_i$ and maximum degree $\Delta_i$, $i\in \{1,2\}$. Then, for every $k \in\{ 1-\delta_1-\delta_2,...,\Delta_1+\Delta_2\}$,
$$\phi_{_k}^d(G_1\square G_2)\ge \alpha(G_1)\alpha(G_2)+\min\{n_1-\alpha(G_1),n_2-\alpha(G_2)\}.$$
\end{remark}

Let $G_1$ be the star graph of order $t+1$ and let $G_2$ be the path graph of order $3$. In this case, $\phi_k^d(G_1\square G_2)=2t+1$ for $k\in \{-1,0\}$. Therefore, the above bound is tight. Even so, Corollary \ref{CoroTh1def(i)} (ii) improves the above bound for the cases where $\phi_{_{k_i}}^d(G_i)>\alpha(G_i)$, for some  $i\in \{1,2\}$.

\begin{theorem}{\em \cite{Yero2013}}\label{th1}
Let $G_i=(V_i,E_i)$ be a simple graph of maximum degree $\Delta_i$,  $i\in \{1,2\}$, and let $S\subseteq  V_1\times V_2 $.
Then the following assertions hold.
\begin{itemize}
\item[{\rm (i)}] If  $P_{V_i}(S)$ is a $k_i$-daf set in $G_i$, then $S$ is a $(k_i+\Delta_j)$-daf set in $G_1\square G_2$, where $j\in \{1,2\}$, $j\ne i.$
\item[{\rm (ii)}]
If for every  $i\in \{1,2\}$, $P_{V_i}(S)$ is a $k_i$-daf set in $G_i$,   then $S$ is a $(k_1+k_2-1)$-daf set in $G_1\square G_2$.

\end{itemize}
\end{theorem}

\begin{corollary}{\em \cite{Yero2013}}\label{CoroTh1def(i)}Let $G_l$ be a graph of order $n_l$, maximum degree
$\Delta_l$ and minimum degree $\delta_l$, with $l\in \{1,2\}$. Then the following assertions hold.
\begin{itemize}
\item[{\rm (i)}]
  For every
$k\in \{\Delta_j-\Delta_i,..., \Delta_i+\Delta_j\}$ $(i,j\in \{1,2\}$, $i\ne j)$,
$$\phi_k^d(G_1\square G_2)\ge
n_j\phi_{k-\Delta_j}^d(G_i).$$

\item[{\rm (ii)}] For every $k_i \in\{ 1-\delta_i,...,\Delta_i\}$, $i\in \{1,2\}$,
$$\phi_{_{k_1+k_2-1}}^d(G_1\square G_2)\ge \phi_{_{k_1}}^d(G_1)\phi_{_{k_2}}^d(G_2)+\min\{n_1-\phi_{_{k_1}}^d(G_1),n_2-\phi_{_{k_2}}^d(G_2)\}.$$
\end{itemize}
\end{corollary}

We emphasize that Corollary \ref{CoroTh1def(i)}  and Proposition \ref{remarktree} lead to infinite families of graphs whose Cartesian product satisfies $\phi_k^d(G_1\square G_2)=n_1n_2.$ For instance, if
$G_1$ is a tree of order $n_1$ and maximum degree $\Delta_1\ge 2$, $G_2$ is a graph of order $n_2$ and  maximum degree $\Delta_2$, and $k\in \{2+\Delta_2,...,\Delta_1+\Delta_2\}$, we have $\phi_k^d(G_1\square G_2)=\phi_{k-\Delta_2}^d(G_1)n_2=n_1n_2.$
In particular, if $G_2$ is a cycle graph, then $\phi_4^d(G_1\square G_2)=n_1n_2.$

Another example of equality in Corollary \ref{CoroTh1def(i)} (ii) is obtained, for instance, taking the
Cartesian product of the star graph $S_t$ of order $t+1$ and the path graph $P_r$ of order
$r$. In that case, for $G_1=S_t$ we have $\delta_1=1$,
$n_1=t+1$ and $\phi_0^d(G_1)=t$, and, for $G_2=P_r$ we have
$\delta_2=1$, $n_2=r$ and $\phi_1^d(G_2)=r-1$. Therefore,
$
\phi_0^d(G_1)\phi_1^d(G_2)+\min\{n_1-\phi_0^d(G_1),n_2-\phi_1^d(G_2)\}=t(r-1)+1.
$
On the other hand, it is not difficult to check that, if we take all
leaves belonging to the copies of $S_t$ corresponding to the first $r-1$
vertices of $G_2$ and we add the vertex of degree $t$ belonging to the
last copy of $S_t$, we obtain a maximum defensive $0$-alliance free set  of cardinality $t(r-1)+1$ in the
graph $G_1\square G_2$, that is,
$\phi_0^d(G_1\square G_2)=t(r-1)+1$. This example also shows
that this  bound is better than the  bound obtained in
Remark \ref{remark1}, which is $t\left\lceil\frac{r}{2}\right\rceil+1$. In this
particular case, both  bounds are equal if and only if $r=2$ or
$r=3$.

\begin{theorem}{\em \cite{Yero2013}}
Let $G_i=(V_i,E_i)$ be a graph and let $S_i\subseteq V_i$, $i\in
\{1,2\}$. If $S_1\times S_2$ is a $k$-daf set in
$G_1\square G_2$ and $S_2$ is a
defensive $k'$-alliance in $G_2$, then $S_1$ is a
$(k-k')$-daf set in $G_1$.
\end{theorem}

Taking into account that $V_2$ is a defensive $\delta_2$-alliance in
$G_2$ we obtain the following result.

\begin{corollary}{\em \cite{Yero2013}}\label{otrocoro}
Let $G_i=(V_i,E_i)$ be a graph,  $i\in
\{1,2\}$. Let  $\delta_2$ be the minimum degree of $G_2$ and let $S_1\subseteq V_1$. If $S_1\times V_2$ is a $k$-daf set in
$G_1\square G_2$, then $S_1$ is a $(k-\delta_2)$-daf set in
$G_1$.
\end{corollary}

By Theorem \ref{th1} (i) and Corollary \ref{otrocoro} we obtain the following result.
\begin{proposition}{\em \cite{Yero2013}}
Let  $G_1$ be a graph of maximum degree $\Delta_1$ and let $G_2$ be a $\delta_2$-regular graph.   For every $k\in \{\delta_2-\Delta_1,...,\Delta_1+\delta_2\}$,  $S_1\times V_2$ is a $k$-daf
set in $G_1\square G_2$ if and only if $S_1$ is a
$(k-\delta_2)$-daf set in $G_1$.
\end{proposition}


\begin{thebibliography}{99}

\bibitem{G-araujo} G. Araujo-Pardo, L. Barri\`ere, Defensive alliances in regular graphs and circulant graphs, \emph{Discrete Applied Mathematics} \textbf{157} (8) (2009) 1924--1931.


\bibitem{chellali-note} M. Bouzefrane, M. Chellali, A note on global alliances in trees, {\em Opuscula Mathematica} \textbf{31} (2) (2011) 153--158.

\bibitem{bouzefrane-glob-def-tree} M. Bouzefrane, M. Chellali, T. W. Haynes, Global defensive alliances in trees, {\em Utilitas Mathematica}
 {\bf 82} (2010) 241--252.

%\bibitem{def-corona} R. G. Eballe, R. M. Aldema, E. M. Paluga, R. F. Rulete, F. P. Jamil, Global defensive alliances in the join, corona and composition of graphs, \emph{Ars Combinatoria}  \textbf{107}  (2012) 225--245.

\bibitem{eroh-bounds} G. Bullington, L. Eroh, S. J. Winters, Bounds concerning the alliance number, {\em Mathematica Bohemica} {\bf 134} (4) (2009) 387--398.

\bibitem{complej1} A. Cami, H. Balakrishnan, N. Deo, R. D. Dutton, On the complexity of finding optimal global alliances, {\it Journal of Combinatorial Mathematics and Combinatorial Computing} {\bf  58} (2006) 23-31.

\bibitem{matamala} R. Carvajal, M. Matamala, I. Rapaport, N. Schabanel, Small alliances in graph, {\em Lectures Notes in Computer Science} {\bf 4708} (2007) 218--227.

\bibitem{Chang2012479}
C. Chang, M. Chia, C. Hsu,  D. Kuo, L. Lai, F. Wang,
Global defensive alliances of trees and Cartesian product of paths and cycles,
 {\it Discrete Applied Mathematics}
 {\bf 160} {4--5} (2012)
  479 -- 487.




\bibitem{chellali} M. Chellali, T. Haynes, Global alliances and independence in trees, {\it Discussiones Mathematicae Graph Theory} {\bf 27} (1) (2007) 19--27.

\bibitem{Chen2013}
 X. Chen  and S. W. Chee, A new upper bound on the global defensive alliance number in trees. \textit{Electronic Journal of Combinatorics}
  \textbf{18} (2011)  Paper 202, 7 pages.

\bibitem{dickson} P. Dickson, K. Weaver, Alliance formation: the relationship between national RD intensity and SME size, {\em Proceedings of ICSB
$50^{th}$ World Conferense D.C.} (2005) 123--154.

\bibitem{enciso-tesis} R. I. Enciso, Alliances in graphs: parameterized algorithms and on partitioning series-parallel graphs. Ph. D. Thesis,
University of Central Florida, Orlando, Florida, USA. 2009.

\bibitem{partitionTrees} L. Eroh, R. Gera,  Global alliance partition in trees, \emph{Journal of Combinatorial Mathematics and Combinatorial Compututing}  \textbf{66} (2008) 161--169.

\bibitem{partitionnumber} L. Eroh, R. Gera, Alliance partition number in graphs, {\em Ars Combinatoria} \textbf{103} (2012) 519--529.

\bibitem{Favaron-ind-dom} O. Favaron, Global alliances and independent domination in some classes of graphs, {\em The Electronic Journal of Combinatorics} {\bf 15} (2008), \#R123.

\bibitem{complej2} H. Fernau, D. Raible, Alliances in graphs: a complexity-theoretic study, Software Seminar SOFSEM 2007, Student Research Forum, Proceedings
Vol. II, 61--70.

\bibitem{flake} G. W. Flake, S. Lawrence, C. L. Giles, Efficient identification of web communities, In {\it Proceedings of the 6th
ACM SIGKDD International Conference on Knowledge Discovery and Data Mining} (KDD-2000) (2000) 150--160.

\bibitem{note} G. H. Fricke, L. M. Lawson, T. W. Haynes, S. M. Hedetniemi, S. T. Hedetniemi, A note on defensive alliances in graphs, {\it Bulletin of the
Institute of Combinatorics and its Applications} {\bf 38} (2003) 37--41.

\bibitem{ararat-ctw} A. Harutyunyan, Some bounds on alliances in trees, \emph{Discrete Applied Mathematics} \textbf{161} (12) (2013) 1739--1746.

\bibitem{GlobalalliancesOne} T. W. Haynes, S. T. Hedetniemi, M. A. Henning, Global defensive alliances in graphs, {\it Electronic Journal of
Combinatorics} {\bf 10} (2003)  139--146.

\bibitem{bookdom1}T. W. Haynes, S. T. Hedetniemi, P. J. Slater, \emph{Fundamentals of Domination in Graphs}, Marcel Dekker, Inc. New York, 1998.



\bibitem{partitionGrid}  T. W. Haynes, J. A. Lachniet,  The alliance partition number of grid graphs, \emph{AKCE International Journal of Graphs
and Combinatorics} \textbf{4} (1) (2007)  51--59.

\bibitem{global-defe-star}  C. J. Hsu, F. H. Wang, Y. L. Wang, Global defensive alliances in star graphs, {\em Discrete Applied Mathematics}  {\bf 157} (8) (2009) 1924--1931.



\bibitem{jamieson-tesis} L. H. Jamieson, Algorithms and complexity for alliances and weighted alliances of varoius types, Ph. D. Thesis, Clemson University, Clemson, SC, USA. 2007.

%\bibitem{jamieson-paralell} L. H. Jamieson, Alliances in generalized series parallel graphs, Proceedings of the Thirty-Ninth Southeastern International Conference on Combinatorics, Graph Theory and Computing, {\em Congressus Numerantium}  {\bf 193} (2008) 157--174.

\bibitem{jamieson-dean} L. H. Jamieson, B. C. Dean, Weighted alliances in graphs, {\em Congressus Numerantium} {\bf 187} (2007) 76--82.

\bibitem{jamieson} L. H. Jamieson, S. T. Hedetniemi, A. A. McRae, The algorithmic complexity of alliances in graphs, \emph{Journal of Combinatorial Mathematics and Combinatorial Computing}  \textbf{68}  (2009), 137-150.


\bibitem{alliancesOne} P. Kristiansen, S. M. Hedetniemi, S. T. Hedetniemi, Alliances in graphs, {\it Journal of Combinatorial
Mathematics and Combinatorial Computing} {\bf 48} (2004) 157--177.

\bibitem{spectral} J. A. Rodr\'{\i}guez, J. M. Sigarreta, Spectral study of alliances in graphs, {\em Discussiones Mathematicae Graph Theory} {\bf 27} (1) (2007) 143--157.

\bibitem{planar} J. A. Rodr\'{\i}guez-Vel\'{a}zquez,  J. M. Sigarreta, Global alliances in planar graphs, {\it  AKCE International Journal of Graphs and
Combinatorics} {\bf 4} (1) (2007) 83--98.

\bibitem{GArs} J. A. Rodr\'{\i}guez-Vel\'{a}zquez,  J. M. Sigarreta, Global defensive $k$-alliances in graphs, \emph{Discrete Applied Mathematics} {\bf  157} (2009)  211--218.

\bibitem{yrs} J. A. Rodr\'{\i}guez-Vel\'{a}zquez, I. G. Yero, J. M. Sigarreta,  Defensive $k$-alliances in graphs, \emph{Applied Mathematics Letters} \textbf{22} (2009)  96--100.

 \bibitem{Rodriguez-Velazquez2011}J. A.  Rodr\'{i}guez-Vel\'azquez,  J. M. Sigarreta, I. G. Yero, S. Bermudo. Alliance free and alliance cover sets.  \textit{Acta Mathematica Sinica, English Series} \textbf{27} (3)  (2011)  497--504.

\bibitem{kdaf1} K. H. Shafique, R. D. Dutton,  Maximum alliance-free and minimum alliance-cover sets, {\it Congressus Numerantium} {\bf 162} (2003) 139--146.

\bibitem{kdaf2} Shafique, K. H., Dutton, R.:  A tight bound on the cardinalities of maximun alliance-free
and minimun alliance-cover sets. {\em J. Combin. Math. Combin. Comput.}, \textbf{56},  139--145 (2006)

\bibitem{kdaf4} K. H. Shafique, Partitioning a Graph in Alliances and its Application to Data Clustering. Ph. D. Thesis, 2004.

\bibitem{kdaf5} K. H. Shafique, R. D. Dutton, On satisfactory partitioning of graphs, {\it Congressus Numerantium} {\bf 154} (2002)  183--194.

\bibitem{tesisiga} J. M. Sigarreta, Alianzas en grafos, Ph. D. Thesis, Universidad Carlos III, Madrid, Spain. 2007.

\bibitem{SBF} J. M. Sigarreta, S. Bermudo, H. Fernau. On the complement graph and defensive $k$-alliances, {\em Discrete Applied Mathematics}
\textbf{157} (8) (2009) 1687--1695.

\bibitem{line} J. M. Sigarreta, J. A. Rodr\'{\i}guez, On defensive alliance and line graphs, {\it Applied Mathematics Letters} {\bf 19} (12) (2006)
1345--1350.

\bibitem{application} G. Szab$\ddot{o}$, T. Cz\'ar\'an, Defensive alliances in spatial models of cyclical population interactions, {\em Physical Reviews E} {\bf 64}, 042902 (2001) 11 pages.




\bibitem{yero2} I. G. Yero, S. Bermudo, J. A. Rodr\'iguez-Vel\'azquez, J. M. Sigarreta, Partitioning a graph into defensive $k$-alliances, \emph{Acta Mathematica Sinica - $($English Series$)$}  \textbf{27} (1) (2011) 73--82.

\bibitem{boundary-def} I. G. Yero, J. A. Rodr\'iguez-Vel\'azquez, Boundary defensive $k$-alliances in graphs, \textit{Discrete Applied Mathematics} \textbf{158} (11) (2010) 1205--1211.

\bibitem{boundary-def-1} I. G. Yero, Contribution to the study of alliances in graphs, Ph. D. Thesis, Universitat Rovira i Virgili, Tarragona, Spain. 2010.

\bibitem{Yero2013}I. G. Yero, J. A. Rodr\'{i}guez-Vel\'{a}zquez and S. Bermudo, Alliance free sets in Cartesian product graphs.  \textit{Discrete Applied Mathematics} \textbf{161} (10--11)  (2013) 1618--1625

\end{thebibliography}
\end{document}